\documentclass[3p,final,number]{elsarticle}
\biboptions{sort&compress}

\usepackage{graphicx}
\graphicspath{{./Figs/}}

\usepackage{amssymb}
\usepackage{amsfonts}
\usepackage{amsmath}

\usepackage{color}

\newtheorem{proposition}{Proposition}

\newcommand{\e}{\mathrm{e}}

\begin{document}
\begin{frontmatter}
\title{On Maximal Total Entropy Production Models for Steady Evaporation of a Calorically Perfect Polyatomic Gas}
\author{Niclas Bernhoff}
\ead{niclas.bernhoff@kau.se}
\author{Eddie Wadbro}
\ead{eddie.wadbro@kau.se}
\address{Department of Mathematics and Computer Science, Karlstad University, Sweden}

\begin{abstract}
  This study investigates the boundary conditions for fluid-dynamic equations at the interface of a gas and its condensed phase during steady evaporation of a polyatomic gas. Evaporation curves illustrating the dependence of the temperature and pressure ratios on the Mach number are considered for a calorically perfect gas whose molecules behave like rigid rotors. At the condensed phase, complete absorption conditions are assumed. Also an extension to cases in which a part of the molecules is diffusely reflected at the condensed phase is also considered. We revisit the half-space evaporation problem for polyatomic gases based on previous results, applied to a slightly modified entropy functional. The structure of the (modified) maximal total entropy production curves is investigated, and simple functions that accurately fit the numerical results are proposed. Functions of the proposed form, with modified coefficients, also fit the evaporation curves previously obtained by different numerical methods surprisingly well. The approximation is performed for different numbers of internal degrees of freedom or ratios of specific heats. Simple functions that depend on the ratio of specific heats are found to fit the obtained evaporation curves for different ratios of specific heats very well.  
\end{abstract}
\end{frontmatter}

\section{Introduction}
The half-space problem of evaporation and condensation for rarefied gases is a classic topic in kinetic theory, with important applications in gas dynamics and phase transition phenomena~\cite{Sone-02, Sone-07}. 
Much work has focused on monatomic gases~\cite{STG-01, BGH-01, BGS-06, BG-21}.
Other studies have extended the analysis to polyatomic gases, accounting for
additional internal degrees of freedom and more realistic physical behavior~%
\cite{FY-06, Fr-07, Be-23d, Be-26, BBW-26}. In a series of seminal papers, K.
Aoki, Y. Sone, and coworkers numerically constructed and theoretically
investigated the evaporation curves and condensation surfaces for monatomic
gases, mainly based on the Bhatnagar--Gross--Krook (BGK) operator~\cite%
{ANSS-91, AS-91, ASY-90, SAY-86, SS-90}. Extensive investigations have
followed, and additional methods have been used to simulate evaporation curves, see the work of Morozov et al.~\cite{MGG-24} and the references therein. Corresponding investigations 
have also been performed for calorically perfect polyatomic gases~\cite{FY-06,
Fr-07}. For monatomic gases, these problems were also analyzed through entropy
inequalities, which provide necessary conditions for the existence of
solutions and allow one to characterize non-equilibrium states without a full
solution of the Boltzmann equation, as shown by Bobylev et al.~\cite%
{BGH-01} and Sone et al.~\cite{STG-01}. Later, these studies were extended to include calorically perfect polyatomic gases~\cite{BBW-26}. In particular, the
principle of maximal entropy production has emerged as a useful tool to
identify physically relevant parameter regimes and compare theoretical
predictions with numerical and experimental data~\cite{BGH-01, BBW-26}.

In this work, we revisit the half-space evaporation problem for calorically
perfect polyatomic gases, based on previous results~\cite{BBW-26}. Our focus
is on the principle of maximal total entropy production and its quantitative
characterization. We investigate the structure of the maximal total entropy
production curves and propose simple functions that accurately fit the
numerical results.

The paper is organized as follows. The kinetic model and any necessary assumptions on it are presented in Section~%
\ref{KM}, while the precise formulation of the half-space problem and
explicit estimates for the macroscopic parameters obtained by using the $%
\mathcal{H}$-theorem~\cite{BBW-26} are addressed in Section~\ref{HSP}. The
(modified) total entropy production and an upper bound for it~\cite{BBW-26}
are presented in Section~\ref{EPE}. In Section~\ref{EM}, exponential models
for the pressure and temperature ratios as functions of the Mach number are
investigated and proposed for fixed ratios of specific heats, while in
Section~\ref{GEM}, a generalized exponential model is obtained with the ratio of
specific heats as a parameter. Mixed boundary conditions, where a fraction of
the particles is diffusely reflected while the rest are completely absorbed
by the condensed phase, are addressed in Section~\ref{SGBC}. Finally, some
concluding remarks are given in Section~\ref{CD}.

\section{Kinetic model\label{KM}}
We consider an ideal gas, consisting of a single species of polyatomic molecules, with
mass $m$, where the polyatomicity is modeled by an internal energy variable $%
I\in $ $\mathbb{R}_{+}$. The distribution functions are nonnegative
functions of the form $f=f\left( t,\boldsymbol{x},\boldsymbol{\xi },I\right) $,
with $t\in \mathbb{R}_{+}$, $\boldsymbol{x}=\left( x,y,z\right) \in \mathbb{R}%
^{3}$, and $\boldsymbol{\xi }=\left( \xi _{1},\xi _{2},\xi _{3}\right) \in 
\mathbb{R}^{3}$. We assume, denoting by $k_{\mathrm{B}}$ the Boltzmann
constant, the (approximative) equation of state $c_{p}-c_{V}=k_{\mathrm{B}%
}/m $ for the specific heat at constant pressure $c_{p}$ and volume $c_{V}$,
respectively. Denoting by $\delta$, with $%
\delta >0$, the number of internal degrees of freedom, the specific
internal energy is given by%
\begin{equation}
e=\frac{3+\delta }{2}\frac{k_{\mathrm{B}}}{m}T\text{.}  \label{sie}
\end{equation}%
Assuming a calorically perfect gas---that is, that the number of internal degrees of freedom does not vary with temperature; hence, $\delta$ is independent of $T$---the specific heat at constant volume is given by%
\begin{equation*}
c_{V}=\frac{\mathrm{d}e}{\mathrm{d}T}=\frac{3+\delta }{2}\frac{k_{\mathrm{B}}%
}{m}\text{,}
\end{equation*}%
and, consequently, the ratio of specific heats is given by 
\begin{equation}
\gamma =\frac{c_{p}}{c_{V}}=\frac{5+\delta }{3+\delta }\text{.}  \label{rsh}
\end{equation}

We consider the real Hilbert space $L^{2}\left( 
\mathrm{d}\boldsymbol{\xi \,}\,\mathrm{d}I\right) $, with
inner product%
\begin{equation*}
\left( f,g\right) =\int_{\mathbb{R}^{3}\times \mathbb{R}_{+}}fg\,\mathrm{d}%
\boldsymbol{\xi \,}\mathrm{d}I\text{ for }f,g\in L^{2}\left( \mathrm{d}%
\boldsymbol{\xi \,}\,\mathrm{d}I\right) \text{.}
\end{equation*}%

The evolution of the distribution functions is (in the absence of external
forces) described by the Boltzmann equation%
\begin{equation}
\frac{\partial f}{\partial t}+\left( \boldsymbol{\xi }\cdot \nabla _{\boldsymbol{%
x}}\right) f=Q_{\delta }\left( f,f\right) \text{,}  \label{BE1}
\end{equation}%
where the collision operator $Q_{\delta }=Q_{\delta }\left( f,f\right) $ is
a quadratic bilinear operator that accounts for changes in the velocities and
internal energies of particles due to binary collisions (assuming that the
gas is rarefied, so that higher-order collisions are negligible).

In this work, a density of states of power-law form $I^{\delta /2-1}$---$%
I^{\delta /2-1}\mathrm{d}I$ representing the number of internal states
between $I$ and $I+\mathrm{d}I$---is considered~\cite{BDLP-94, GP-23, Be-23b}
to recover the proper form~\eqref{sie} of the specific internal energy~\cite%
{AMR-24, BBCG-26}. From quantum mechanical results~\cite{Anderson-03,
Anderson-06, Atkins-10}, when the molecules are modeled as rigid rotors, the density of states of
power-law form $I^{\delta /2-1}$ can be motivated in the case of rotational energy~\cite%
{BBCG-26}, where, for example, $\delta =2$ for linear molecules and $\delta =3$ for
spherical tops (at least approximately). For a calorically perfect gas, as
the number of internal degrees of freedom $\delta $ is constant, the density
of states of power-law form $I^{\delta /2-1}$ may, at least for some
purposes, be physically relevant~\cite{DPT-21, BBCG-26}. 
For thermally perfect gases, the number of internal degrees of freedom varies with the temperature, and other densities of states may have to be considered to capture those gases. 
Another approach to capturing thermally perfect gases is to introduce additional discrete or continuous variable(s) for the vibrational part of the internal energy. 
For some purposes, it may also be satisfactory to apply the density of states of power-law form $I^{\delta/2-1}$ for thermally perfect gases, but with an average value of the number of internal degrees of freedom $\delta >0$.

\subsection{Microscopic model}

A collision (localized in space and time) can be represented by the
microscopic velocities and internal energies of the colliding molecules
before and after the collision, denoted by $\left( \boldsymbol{\xi }%
,I\right) $ and $\left( \boldsymbol{\xi }_{\ast },I_{\ast }\right) $, and $%
\left( \boldsymbol{\xi }^{\prime },I^{\prime }\right) $ and $\left( 
\boldsymbol{\xi }_{\ast }^{\prime },I_{\ast }^{\prime }\right) $,
respectively. The momentum and total energy conservation of
the collision read 
\begin{align*}
\boldsymbol{\xi }+\boldsymbol{\xi }_{\ast }& =\boldsymbol{\xi }^{\prime }+%
\boldsymbol{\xi }_{\ast }^{\prime }\text{,} \\
\frac{m}{2}\left\vert \boldsymbol{\xi }\right\vert ^{2}+\frac{m}{2}%
\left\vert \boldsymbol{\xi }_{\ast }\right\vert ^{2}+I+I_{\ast }& =\frac{m}{2%
}\left\vert \boldsymbol{\xi }^{\prime }\right\vert ^{2}+\frac{m}{2}%
\left\vert \boldsymbol{\xi }_{\ast }^{\prime }\right\vert ^{2}+I^{\prime
}+I_{\ast }^{\prime }\text{.}
\end{align*}%
In the center-of-mass frame, energy conservation reads%
\begin{equation*}
E:=\frac{m}{4}\left\vert \boldsymbol{\xi }-\boldsymbol{\xi }_{\ast
}\right\vert ^{2}+I+I_{\ast }=\frac{m}{4}\left\vert \boldsymbol{\xi }%
^{\prime }-\boldsymbol{\xi }_{\ast }^{\prime }\right\vert ^{2}+I^{\prime
}+I_{\ast }^{\prime }=:E^{\prime }\text{,}
\end{equation*}%
defining the total energy $E$ in the center-of-mass frame.
For resonant collisions, see~\cite{BRS-24, AB-26, Be-26}, the total energy conservation splits into kinetic and internal energy conservation laws%
\begin{equation*}
\left\vert \boldsymbol{\xi }\right\vert ^{2}+\left\vert \boldsymbol{\xi }%
_{\ast }\right\vert ^{2}=\left\vert \boldsymbol{\xi }^{\prime }\right\vert
^{2}+\left\vert \boldsymbol{\xi }_{\ast }^{\prime }\right\vert ^{2}\text{
and }I+I_{\ast }=I^{\prime }+I_{\ast }^{\prime }\text{,}
\end{equation*}%
or equivalently, conservation of relative velocity and vanishing internal energy gap%
\begin{equation*}
\left\vert \boldsymbol{\xi }-\boldsymbol{\xi }_{\ast }\right\vert
=\left\vert \boldsymbol{\xi }^{\prime }-\boldsymbol{\xi }_{\ast }^{\prime
}\right\vert \text{ and }I^{\prime }+I_{\ast }^{\prime }-I-I_{\ast }=0\text{.%
}
\end{equation*}

\subsection{Macroscopic quantities and compressible Euler system}

Macroscopic quantities---that is, the number density of molecules $n$, the mass
density $\rho $, the flow velocity $\boldsymbol{u}=\left(
u_{1},u_{2},u_{3}\right) $, the temperature $T$, and the pressure $p$---are
defined by%
\begin{align*}
n& =\left( 1,f\right) \text{, }\rho =mn=\left( m,f\right) \text{, }u_{i}=%
\frac{1}{n}\left( \xi _{i},f\right) \text{, }i=1,2,3\text{, }\\
T&=\frac{2}{\left( 3+\delta \right) nk_{\mathrm{B}}}\left( \frac{m}{2}%
\left\vert \boldsymbol{\xi }-\boldsymbol{u}\right\vert ^{2}+I,f\right) \text{,
and } p=nk_{\mathrm{B}}T=\frac{2}{3+\delta }\left( \frac{m}{2}\left\vert 
\boldsymbol{\xi }-\boldsymbol{u}\right\vert ^{2}+I,f\right) \text{.}
\end{align*}

In the hydrodynamic limit, that is, when the Knudsen number tends to zero, the
overall evolution of the macroscopic quantities is governed by the
compressible Euler equations (in the absence of external forces)%
\begin{align*}
\frac{\partial \rho }{\partial t}+\nabla _{\boldsymbol{x}}\cdot \left( \rho 
\boldsymbol{u}\right) &= 0\text{,} \\
\rho \frac{\partial }{\partial t}\boldsymbol{u}+\rho \left( \boldsymbol{u}\cdot
\nabla _{\boldsymbol{x}}\right) \boldsymbol{u}+\frac{k_{\mathrm{B}}}{m}\nabla _{%
\boldsymbol{x}}\left( \rho T\right) &= 0\text{,} \\
\frac{\partial T}{\partial t}+\boldsymbol{u}\cdot \nabla _{\boldsymbol{x}}T+\left(
\gamma -1\right) T\,\nabla _{\boldsymbol{x}}\cdot \boldsymbol{u} &=0\text{,}
\end{align*}%
which can be obtained through a Chapman--Enskog expansion of the Boltzmann
equation~\eqref{BE1} for the density of states of power-law form $I^{\delta
/2-1}$~\cite{BBBD-18}.

The characteristics of the corresponding one-dimensional Euler system are $\left\{ u-c,u,u+c\right\} $, where%
\begin{equation*}
c=\sqrt{\frac{\gamma p}{\rho }}=\sqrt{\frac{\gamma k_{\mathrm{B}}T}{m}}
\end{equation*}%
denotes the speed of sound.

\subsection{Assumed properties of the collision operator\label{MPCO}}

In this section, we state the assumed properties of the collision operator $%
Q_{\delta }\left( f,f\right) $.

There are five conservation laws~\cite{BDLP-94}%
\begin{equation}
\left( \psi \left( \boldsymbol{\xi },I\right) ,Q_{\delta
}(f,f)\right) =0\text{ for }\psi \left( \boldsymbol{\xi }, I\right)
\in \left\{ 1,\xi _{1},\xi _{2},\xi _{3},m\left\vert \boldsymbol{\xi }%
\right\vert ^{2}+2I\right\} \text{.}  \label{conslaw}
\end{equation}%
The collision operator satisfies the $\mathcal{H}$-theorem~\cite{BDLP-94, Be-23b}, which states that%
\begin{equation}
\left( \log \left( I^{1-\delta /2}f\right) ,Q_{\delta }(f,f)\right) \leq 0%
\text{,}  \label{HT}
\end{equation}%
where equality holds in inequality~\eqref{HT} if and only if 
\begin{equation*}
Q_{\delta }(f,f)=0\text{,}
\end{equation*}%
or, if and only if there exist $n\geq 0$, $\boldsymbol{u}\in \mathbb{R}^{3}$,
and $T>0$, such that for almost every $\left( \boldsymbol{\xi },I\right) \in 
\mathbb{R}^{3}\times \mathbb{R}_{+}$%
\begin{align*}
f =M\left( \boldsymbol{\xi },I\right) =M_{\mathrm{tr}}\left( \boldsymbol{%
\xi }\right) M_{\mathrm{int}}\left( I\right) \text{, with }M_{\mathrm{tr}}&=%
\frac{nm^{3/2}}{\left( 2\pi k_{\mathrm{B}}T\right) ^{3/2}}\exp \left( -\frac{%
m|\boldsymbol{\xi }-\boldsymbol{u}|^{2}}{2k_{\mathrm{B}}T}\right) 
\\ \text{and }M_{\mathrm{int}}&=\frac{I^{\delta /2-1}}{\Gamma \left( \delta
/2\right) \left( k_{\mathrm{B}}T\right) ^{\delta /2}}\exp \left( -\frac{I}{%
k_{\mathrm{B}}T}\right) \text{,}
\end{align*}%
is a Maxwellian distribution.
In the expression above, $\Gamma$ is the usual Gamma function.

\section{Half-space problem of evaporation and condensation \label{HSP}}

We consider the stationary Boltzmann equation in one spatial dimension; that
is, $f$ depends only on a single space variable, henceforth denoted by $x>0$, while retaining dependence on the three velocity variables $\boldsymbol{\xi }=\left( \xi _{1},\xi_{2},\xi
_{3}\right) $. Then 
\begin{equation}
\xi _{1}\dfrac{\partial f}{\partial x}=Q_{\delta }(f,f),  \label{HSP1}
\end{equation}%
where $f=f(x,\boldsymbol{\xi },I)$ represents the distribution function of the
molecules at position $x\in \mathbb{R}_{+}$, with velocity $\boldsymbol{\xi }%
\in \mathbb{R}^{3}$ and internal energy $I\in \mathbb{R}_{+}$, and $\delta
>0 $ denotes the number of internal degrees of freedom.

\subsection{Boundary conditions}

Introduce the notation (where $f=f(\boldsymbol{\xi})$ may depend on
more variables than $\boldsymbol{\xi}\in \mathbb{R}^{3}$) 
\begin{equation*}
f_{\pm }(\boldsymbol{\xi })=f_{\pm }(\xi _{1},\xi _{2},\xi _{3})=f(\pm \xi
_{1},\xi _{2},\xi _{3})\hspace*{3mm}\text{for }\xi _{1}>0\text{.}
\end{equation*}%
Assuming complete absorption with a non-drifting incoming Maxwellian
distribution $M_{0+}$ from the condensed phase and an equilibrium distribution $M_{\infty }$ being approached at the far end, we obtain the boundary conditions%
\begin{equation}
f_{+}(0,\boldsymbol{\xi}, I)=M_{0+}\text{ and } f(x,\boldsymbol{\xi }, I)\rightarrow M_{\infty }\text{ as $x\rightarrow \infty$,}
\label{BC}
\end{equation}%
with
\begin{align*}
M_{0} & =
M\vert_{\left(n,\boldsymbol{u},T\right) = \left( n_{0},\boldsymbol{0},T_{0}\right)} = 
\left. M_{\mathrm{tr}}\right\vert_{\left( n,\boldsymbol{u},T\right) =\left(n_{0},\boldsymbol{0},T_{0}\right) }
\left. M_{\mathrm{int}}\right\vert _{T=T_{0}} 
\text{ and} \\
M_{\infty }& =
M\vert _{\left( n,\boldsymbol{u},T\right) =\left(
n_{\infty },\boldsymbol{u}_{\infty },T_{\infty }\right) }=
\left. M_{\mathrm{tr}}\right\vert _{\left( n,\boldsymbol{u},T\right) =\left( n_{\infty },\boldsymbol{u}_{\infty },T_{\infty }\right) }
\left. M_{\mathrm{int}}\right\vert_{T=T_{\infty }}\text{,}
\end{align*}%
where we assume that $\boldsymbol{u}_{\infty }=\left(u,0,0\right) $.
The Mach number at the far end is defined as%
\begin{equation*}
\boldsymbol{\mathcal{M}}_{\infty} = 
\dfrac{\boldsymbol{u}_{\infty }}{c} = 
\left( \mathcal{M},0,0\right)
\text{, with }
\mathcal{M}=\dfrac{u}{c}=\sqrt{\dfrac{m}{\gamma k_{\mathrm{B}}T_{\infty }}}u\text{.}
\end{equation*}%
In this paper, we limit our attention to subsonic evaporation $0<u\leq c$
or, equivalently, $0<\mathcal{M}\leq 1$. We introduce the relative pressure
and relative temperature%
\begin{equation*}
\overline{p}=\frac{p_{\infty }}{p_{0}}\text{ and }\overline{T}=\frac{%
T_{\infty }}{T_{0}}\text{.}
\end{equation*}%
Furthermore, we drop the bars and denote by $p$ and $T$ the relative
pressure $\overline{p}$ and relative temperature $\overline{T}$, respectively.

In Section~\ref{SGBC}, we consider generalized boundary conditions~\cite%
{Fr-07, STG-01} where a fraction (but not all) of the molecules are
diffusely re-emitted into the gas after interaction with the interface
through
\begin{equation}
f_{+}(0,\boldsymbol{\xi },I)=\left( \sigma _{e}+(1-\sigma _{e})%
\frac{N_{0}}{n_{0}}\right) M_{0}\text{, where }N_{0}=-\sqrt{\frac{2\pi m}{k_{%
\mathrm{B}}T_{0}}}\int_{\mathbb{R}_{-}^{3}}\xi _{1}f\,\mathrm{d}\boldsymbol{%
\xi }\text{,}  \label{BC1}
\end{equation}%
for $0<\sigma _{e}\leq 1$, with $M_{0}=\left. M\right\vert _{\left( n,%
\boldsymbol{u},T\right) =\left( n_{0},\boldsymbol{0},T_{0}\right) }$. Here $\mathbb{R%
}_{-}^{3}=\left\{ \left. \boldsymbol{\xi }\in \mathbb{R}^{3}\right\vert \xi
_{1}<0\right\} $.

\subsection{Necessary conditions on the boundary data \label{NCBD}}

Applying the conservation laws~\eqref{conslaw}, we obtain%
\begin{align*}
(1,\xi _{1}f) &=n_{\infty }u\text{, }(\xi _{1},\xi _{1}f)=n_{\infty }\left(
u^{2}+k_{\mathrm{B}}T_{\infty }/m\right) \text{, }\\(\xi _{2},\xi _{1}f)&=(\xi
_{3},\xi _{1}f)=0\text{, and } 
(\left\vert \boldsymbol{\xi }\right\vert ^{2}+2I/m,\xi
_{1}f)=n_{\infty }u\left( u^{2}+\left( 5+\delta \right) k_{\mathrm{B}%
}T_{\infty }/m\right) \text{.}
\end{align*}%
Moreover, we introduce the $\mathcal{H}$-functional 
\begin{equation}
\Psi \left( f\right) =\left( \log \left( I^{1-\delta /2}f\right) ,\xi
_{1}f\right) \text{,}  \label{HF}
\end{equation}%
and, by the $\mathcal{H}$-theorem~\eqref{HT}, obtain
\begin{equation}
\frac{\mathrm{d}\Psi \left( f\right) }{\mathrm{d}x}=\left( \log \left( I^{1-\delta /2}f\right)
,Q_{\delta }(f,f)\right) \leq 0\text{.}  \label{MHT}
\end{equation}%
Hence, for $u>0$ or, equivalently, $\mathcal{M}>0$, we have that%
\begin{equation}
\frac{\mathrm{d}\widetilde{\Psi }\left( f\right) }{\mathrm{d}x}\geq 0\text{ for }\widetilde{%
\Psi }\left( f\right) =-\frac{\Psi \left( f\right) }{(1,\xi _{1}f)}=-\frac{%
\left( \log \left( I^{1-\delta /2}f\right) ,\xi _{1}f\right) }{n_{\infty }u}%
\text{.}  \label{NE}
\end{equation}%
Note that 
\begin{equation*}
n_{\infty }u=p_{\infty }\sqrt{\dfrac{\gamma }{mk_{\mathrm{B}}T_{\infty }}}%
\mathcal{M}\text{.}
\end{equation*}

Introducing the ``half''-moments $N_{i}^{\pm }$, $i=1,\dotsc ,5$, by 
\begin{equation*}
\int\limits_{\mathbb{R}_{+}^{3}\times \mathbb{R}_{+}}\xi _{1}%
\begin{Bmatrix}
1 \\ 
\xi _{1} \\ 
\xi _{2} \\ 
\xi _{3} \\ 
\left\vert \boldsymbol{\xi }\right\vert ^{2}+2I/m%
\end{Bmatrix}%
f_{\pm }\left( 0,\boldsymbol{\xi },I\right) \,\mathrm{d}\boldsymbol{\xi }\,%
\mathrm{d}I=:%
\begin{Bmatrix}
N_{1}^{\pm } \\ 
N_{2}^{\pm } \\ 
N_{3}^{\pm } \\ 
N_{4}^{\pm } \\ 
N_{5}^{\pm }%
\end{Bmatrix}%
\text{,}
\end{equation*}%
where $\mathbb{R}_{+}^{3}=\left\{ \left. \boldsymbol{\xi }\in \mathbb{R}%
^{3}\right\vert \xi _{1}>0\right\} $, we obtain%
\begin{equation}
\begin{aligned} N_{1}^{-} &=N_{1}^{+}-(1,\xi
_{1}f)=\frac{p_{0}}{\sqrt{mk_{\mathrm{B}}T_{0}}}\left( \frac{1}{\sqrt{2\pi
}}-p\sqrt{\dfrac{\gamma }{T}}\mathcal{M}\right) \geq 0\text{,} \\ N_{2}^{-} &=(\xi
_{1},\xi _{1}f)-N_{2}^{+}=\frac{p_{0}}{m}\left( p\left( 1+\gamma
\mathcal{M}^{2}\right) -\frac{1}{2}\right) \geq 0\text{, and} \\
N_{5}^{-}&=N_{5}^{+}-(\left\vert \boldsymbol{\xi }\right\vert ^{2}+2I/m,\xi
_{1}f) \\ &=\frac{p_{0}\sqrt{k_{\mathrm{B}}T_{0}}}{m^{3/2}}\left(
\frac{4+\delta }{\sqrt{2\pi }}-p\sqrt{\gamma T}\mathcal{M}\left( 5+\delta
+\gamma \mathcal{M}^{2}\right) \right) \geq 0\text{.} \end{aligned}
\label{NMR}
\end{equation}

Note that $N_{3}^{\pm }=N_{4}^{\pm }=0$, and for simplicity, we assume that the flows are symmetric such that
\begin{equation*}
f=f\left( x,\boldsymbol{\xi },I\right) =f\left( x,\xi _{1},r,I\right)
\text{ for }
r=\sqrt{\xi _{2}^{2}+\xi _{3}^{2}}\text{.}
\end{equation*}
The following proposition summarizes the results obtained in our previous work~\cite{BBW-26} on the
necessary conditions on the boundary data for the existence of solutions to
the half-space problem for evaporation.

\begin{proposition}
\label{P1}For the half-space problem~\eqref{HSP1}, with boundary conditions~%
\eqref{BC}, to admit a solution, the following relations between
the parameters of the two Maxwellians at the condensed interface and the
uniform phase at infinity are forced to be fulfilled:
\begin{itemize}
\item[i)] for all $u>0$ or equivalently for all $\mathcal{M}_{\infty }>0$
\begin{equation}
\begin{aligned} 
p\mathcal{M} &\leq \sqrt{\frac{T}{2\pi \gamma }}\text{,} \\
p\mathcal{M} \sqrt{T} &\leq \frac{%
1 }{\gamma \sqrt{2\pi \gamma }}\frac{4+\delta }{3+\delta +\mathcal{M}^{2}}\text{,} \\
\frac{1}{2\left( 1+\gamma \mathcal{M}^{2}\right) } &\leq p\leq \frac{1}{%
\left( 1+\dfrac{\mathcal{M}^{2}}{3+\delta }\right) ^{\left( 5+\delta \right)
/2}}\text{;}
\end{aligned}  \label{NCP1}
\end{equation}
\item[ii)] for  $u=0$ or equivalently $\mathcal{M}_{\infty }=0$
\begin{equation}
p=T=1\text{.}  \label{TC}
\end{equation}%
Moreover, under the condition~\eqref{TC} there exists a unique solution 
\begin{equation*}
  f=f(x,\boldsymbol{\xi },I) =  M_{0}(\boldsymbol{\xi},I)\text{.}
\end{equation*}
\end{itemize}
\end{proposition}

\section{Total entropy production estimate \label{EPE}}

We denote by the functional%
\begin{equation}
\widetilde{\mathfrak{D}}(f)=-\frac{1}{(1,\xi _{1}f)}\int_{0}^{\infty }\left(
\log \left( I^{1-\delta /2}f\right) ,Q_{\delta }(f,f)\right) \,\mathrm{d}%
x\geq 0  \label{MEP}
\end{equation}%
the total entropy production. Note that this definition differs from that 
of the total entropy production $\mathfrak{D}(f)$ in our previous work~\cite{BBW-26} (as well as in the work by Bobylev et al.~\cite{BGH-01}) by a
factor (independent of $x$) $\left( n_{\infty }u\right) ^{-1}$, that is, $\widetilde{\mathfrak{D}}=\mathfrak{D}(f)/\left( n_{\infty }u\right)$, in view of
the new ``entropy functional'' $\widetilde{\Psi }\left( f\right) $ defined in~\eqref{NE}.

Then, by direct implementation of our previous results~\cite{BBW-26}, we have an
upper bound for the total entropy production%
\begin{equation*}
0\leq \widetilde{\mathfrak{D}}(f)\leq \frac{\sqrt{T}}{\sqrt{\gamma }p\mathcal{M}}\Lambda
(p,T,\mathcal{M})\text{,}
\end{equation*}%
where%
\begin{multline*}
\Lambda (p,T,\mathcal{M})=\frac{1}{\sqrt{2\pi }}\log \frac{T^{\left(
5+\delta \right) /2}}{p\sqrt{e}} \\
+\left( \sqrt{\frac{\gamma }{T}}p\mathcal{M-}\frac{1}{\sqrt{2\pi }}\right)
\log \left( \frac{2\left( \sqrt{\pi }\left( p\left( 1+\gamma \mathcal{M}%
^{2}\right) -1/2\right) \right) ^{6+\delta }T^{\left( 5+\delta \right)
/2}\Delta e^{-\widetilde{\theta }(s)}}{p\left( 4+\delta -p\sqrt{2\pi \gamma T%
}\mathcal{M}\left( 5+\delta +\gamma \mathcal{M}^{2}\right) \right)
^{5+\delta }}\right) \geq 0\text{,}
\end{multline*}%
with%
\begin{align*}
\widetilde{\theta }(s)& =\frac{1}{2}-s^{2}-s\frac{I_{0}(s)}{I_{1}(s)}=\frac{1%
}{2}-s\frac{I_{2}(s)}{I_{1}(s)}=\frac{1}{2}-\frac{s^{2}+s\sqrt{s^{2}+2\left(
4+\delta \right) \Upsilon }}{2\Upsilon }\text{,} \\
\Delta & =\left( \Upsilon +s^{2}\left( 2\Upsilon -1\right) -s\sqrt{%
s^{2}+2\left( 4+\delta \right) \Upsilon }\right) \left( s+\sqrt{%
s^{2}+2\left( 4+\delta \right) \Upsilon }\right) ^{4+\delta }\text{, and} \\
\Upsilon & =\frac{2\left( 1-p\sqrt{2\pi \gamma /T}\mathcal{M}\right) \left(
4+\delta -p\sqrt{2\pi \gamma T}\mathcal{M}\left( 5+\delta +\gamma \mathcal{M}%
^{2}\right) \right) }{\pi \left( 2p\left( 1+\gamma \mathcal{M}^{2}\right)
-1\right) ^{2}}>1\text{.}
\end{align*}%
The parameter $s\in (-\infty ,\infty )$ can be obtained by solving the
equation 
\begin{equation*}
\frac{N_{1}^{-}N_{5}^{-}}{\left( N_{2}^{-}\right) ^{2}}=\frac{I_{1}(s)\left[
I_{3}(s)+I_{1}(s)\left( 1+\delta /2\right) \right] }{I_{2}^{2}(s)}\text{,}
\end{equation*}%
in which 
\begin{equation*}
I_{n}(s)=\int\limits_{0}^{\infty }z^{n}e^{-(z-s)^{2}}\, \mathrm{d}z=\int\limits_{-s}^{\infty }\left( z+s\right) ^{n}e^{-z^{2}}\,\mathrm{d}z%
\text{ for }n\in \left\{ 0,1,...\right\} \text{,}
\end{equation*}%
and the quantities $N_{1}^{-}$, $N_{2}^{-}$, and $N_{5}^{-}$ are given by relations~\eqref{NMR}.

The allowed physical domain of positive total entropy production in the $\left(
p,T,\mathcal{M}\right) $-space is bounded by the surface 
\begin{equation*}
S:\widetilde{\Lambda }(p,T,\mathcal{M})=0\text{,}
\end{equation*}%
where%
\begin{equation*}
\widetilde{\Lambda }(p,T,\mathcal{M})=\frac{\sqrt{T}}{\sqrt{\gamma }p%
\mathcal{M}}\Lambda (p,T,\mathcal{M})\text{.}
\end{equation*}%
For any fixed $\delta $ and $\mathcal{M}$, we define the pair $\left(
p_{\#}^{\delta }\left( \mathcal{M}\right) ,T_{\#}^{\delta }\left( \mathcal{M}%
\right) \right) $ as the relative pressure and temperature values that maximize the
total entropy production~\eqref{MEP} by
\begin{equation*}
\widetilde{\Lambda }\left( p_{\#}^{\delta }\left( \mathcal{M}\right)
,T_{\#}^{\delta }\left( \mathcal{M}\right) ,\mathcal{M}\right) =\max_{p,T}%
\widetilde{\Lambda }\left( p,T,\mathcal{M}\right) .
\end{equation*}

\section{Exponential model fitting\label{EM}}
In this section, our aim is to find a simple function that fits the maximal total entropy production curves for different numbers of internal degrees of freedom $\delta$ and Mach numbers $\mathcal{M}$.

To achieve this, we adopt a regression approach with the goal of identifying functional forms that capture the observed trends.
Our methodology comprises three steps: (i) sampling maximal total entropy production points $(p^\delta_\#(\mathcal{M}), T^\delta_\#(\mathcal{M}))$ for $\mathcal{M} \in (0,1]$ using numerical optimization; (ii) proposing candidate functional forms; and (iii) estimating parameters using nonlinear least-squares fitting and evaluating the accuracy using $R^2$ and mean absolute percentage error (MAPE).

\subsection{Model 1: Basic exponential model}
The first models we fit are based on the observation that the maximal total entropy curves exhibit a shape reminiscent of exponential decay when plotted against the Mach number $\mathcal{M}$.
Based on this observation, we make the exponential ansatz that these curves can be approximated by
\begin{equation}
  p_\dag^\delta(\mathcal{M}) = \e^{-\beta \mathcal{M}} \quad\text{and}\quad
  T_\dag^\delta(\mathcal{M}) = \e^{-\alpha \mathcal{M}},
  \label{eq:fit1}
\end{equation}
where $\alpha>0$ and $\beta>0$ are the model parameters.
For each $\delta \in \{0,2,3,5\}$, we individually fit the coefficients $\alpha$ and $\beta$ in the least-squares sense to the points of the maximal total entropy production for $\mathcal{M}\in\left\{0.01,0.02,\ldots,1\right\}$.
We denote the optimized parameters in the least-squares sense by $\alpha^\delta_\dag$ and $\beta^\delta_\dag$. 
Table~\ref{tab:fit1} summarizes the fitted coefficients and the model performance, evaluated using the coefficient of determination ($R^2$) and the mean absolute percentage error (MAPE).
\begin{table}
\centering
\begin{tabular}{|c|c|c|c|c|c|c|}
\hline
 $\delta$ & $\alpha^\delta_\dag$ & $\beta^\delta_\dag$ & $R^2$ for $T$ & $R^2$ for $p$ & MAPE for $T$ & MAPE for $p$ \\
\hline
0 & 0.4344 & 1.7093 &  0.9992 &  0.9948 &  0.3297 &  3.7658 \\
2 & 0.2735 & 1.5502 &  0.9982 &  0.9954 &  0.3044 &  3.1029 \\
3 & 0.2309 & 1.5079 &  0.9979 &  0.9956 &  0.2797 &  2.9429 \\
5 & 0.1760 & 1.4535 &  0.9974 &  0.9958 &  0.2362 &  2.7476 \\
\hline
\end{tabular}
\caption{Summary of model quality for the basic exponential model~\eqref{eq:fit1}.}
\label{tab:fit1}
\end{table}

\begin{figure}
\centering
\includegraphics[scale=1]{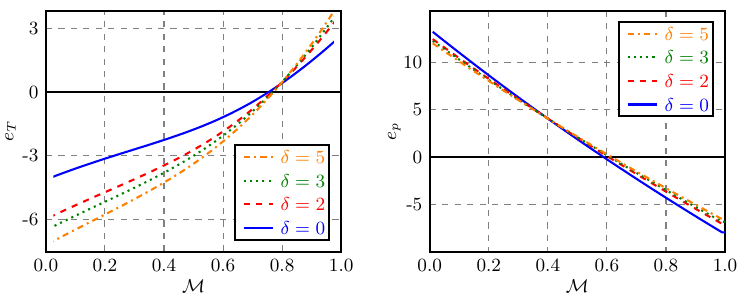}
\caption{Required corrections for a perfect fit for the basic exponential fit for each $\delta \in \{0,2,3,5\}$.}
\label{fig:fit1}
\end{figure}

Although this basic exponential model achieves excellent accuracy for $T$ and very good accuracy for $p$, systematic deviations remain.
To quantify these deviations, we define the correction functions $e_p=e_p(\mathcal{M})$ and $e_T=e_T(\mathcal{M})$, such that 
\begin{equation}
  p_\#^\delta(\mathcal{M}) = \e^{-\beta^\delta_\dag(1+e_p) \mathcal{M}} \quad\text{and}\quad
  T_\#^\delta(\mathcal{M}) = \e^{-\alpha^\delta_\dag(1+e_T) \mathcal{M}}.
  \label{eq:fit1corr}
\end{equation}
Figure~\ref{fig:fit1} illustrates these correction functions.
This suggests that although the fit is very good in terms of $R^2$ and MAPE, there is a systematic error.
For $p$ the required correction appears to be linearly dependent on $\mathcal{M}$, whereas for $T$ the required correction is non-linear---something that we will use as a basis for improving our models. 

\subsection{Model 2: Improved exponential model}
\begin{table}
\centering
\begin{tabular}{|c|c|c|c|c|}
\hline
 $\delta$ & $\alpha^\delta_\ddag$ & $\beta^\delta_\ddag$ & MAPE for $T$ & MAPE for $p$ \\
 \hline
0 & (0.4201, 0.0255) & (1.9241,  0.3599) &  0.0115 &  0.1262 \\
2 & (0.2601, 0.0234) & (1.7345,  0.3007) &  0.0088 &  0.0940 \\
3 & (0.2184, 0.0214) & (1.6846,  0.2862) &  0.0079 &  0.0866 \\
5 & (0.1654, 0.0180) & (1.6206,  0.2683) &  0.0065 &  0.0775 \\
\hline
\multicolumn{5}{|c|}{BBW MTEP from~\cite{BBW-26}} \\
 \hline 
 0 & (0.4267, 0.0216) & (1.9231,  0.3597) &  0.0290 &  0.1281 \\
 2 & (0.2660, 0.0202) & (1.7332,  0.3003) &  0.0212 &  0.0968 \\
 3 & (0.2239, 0.0186) & (1.6832,  0.2858) &  0.0185 &  0.0896 \\
 5 & (0.1700, 0.0157) & (1.6191,  0.2679) &  0.0149 &  0.0810 \\
 \hline

\end{tabular}
\caption{Summary of model quality for the improved exponential model for the modified maximal total entropy production considered in this paper above and the maximal total entropy production considered by Bernhoff, Brull, and Wadbro in~\cite{BBW-26} below.}
\label{tab:fit2}
\end{table}
Based on the observation that there is a systematic deviation when using the basic exponential fit and on the shapes of the required corrections displayed in Figure~\ref{fig:fit1}, we introduce an improved exponential model that incorporates higher-order terms in $\mathcal{M}$ to capture nonlinear effects.
That is, we make the extended exponential ansatz
\begin{equation}
  p_\ddag^\delta(\mathcal{M}) = \e^{-\beta_1 \mathcal{M} + \beta_2 \mathcal{M}^2} \quad\text{and}\quad
  T_\ddag^\delta(\mathcal{M}) = \e^{-\alpha_1 \mathcal{M} - \alpha_2 \mathcal{M}^3},
  \label{eq:fit2}
\end{equation}
where $\alpha_1$, $\alpha_2$, $\beta_1$, and $\beta_2$ are the model parameters.

Then $(p,T)=(p_\ddag^\delta(\mathcal{M}),T_\ddag^\delta(\mathcal{M}))$ is the solution to the system
\begin{equation*}
  \frac{\mathrm{d}p}{p} = -\left(\beta_1  - 2\beta_2 \mathcal{M}\right)\mathrm{d}\mathcal{M}\text{, }\quad
  \frac{\mathrm{d}T}{T} = -\left(\alpha_1 + 3\alpha_2 \mathcal{M}^2\right)\mathrm{d}\mathcal{M},\quad p(0)=T(0)=1\text{.}
\end{equation*}
As before, we fit the parameters individually for each $\delta \in \{0,2,3,5\}$ to the maximal total entropy production curves sampled at $\mathcal{M}\in\left\{0.01,0.02,\ldots,1\right\}$ in the least-squares sense.
We denote the optimized parameters in the least-squares sense by $\alpha^\delta_{\ddag}$ and $\beta^\delta_\ddag$. 
Table~\ref{tab:fit2} reports the fitted coefficients and the MAPE values obtained for the modified maximum total entropy production considered in this paper (upper table), but also the one considered by Bernhoff, Brull, and Wadbro (BBW)~\cite{BBW-26} (lower table); the $R^2$-values are omitted from the table because they round to 1.0000 for all cases, indicating near-perfect fits.

\subsection{Exponential model applied on other data sets}
\begin{table}
\centering
\begin{tabular}{|c|c|c|c|c|c|c|}
\hline
 & $\delta$ & $n$ & $\alpha^\delta_\ddag$ & $\beta^\delta_\ddag$ & MAPE $T$ & MAPE $p$ \\
 \hline 
MTEP & 0 & 100 & (0.4201, 0.0255) & (1.9241,  0.3599) &  0.0115 &  0.1262 \\
DSMC & 0   &  12 & (0.4126, 0.0282) & (1.9350,  0.3650) &  0.0238 &  0.1633 \\ 
 BGK & 0   &  12 & (0.4102, 0.0308) & (1.9281,  0.3592) &  0.0113 &  0.1483 \\ 
S-model & 0   &  12 & (0.4039, 0.0334) & (1.9289,  0.3597) &  0.0122 &  0.1487 \\
 MM & 0   &  12 & (0.4045, -0.0027) & (1.9219,  0.3470) &  0.0034 &  0.1690 \\
MTEP & 2 & 100 & (0.2601, 0.0234) & (1.7345,  0.3007) &  0.0088 &  0.0940 \\
DSMC & 2   &  11 & (0.2604, 0.0206) & (1.7571,  0.3213) &  0.0774 &  0.1522 \\
Holway & 2   &   9 & (0.2572, 0.0243) & (1.7459,  0.3074) &  0.0180 &  0.1216 \\
 MM & 2   & 201 & (0.2472, -0.0006) & (1.7343,  0.2906) &  0.0000 &  0.1015 \\ 
MTEP & 3 & 100 & (0.2184, 0.0214) & (1.6846,  0.2862) &  0.0079 &  0.0866 \\
DSMC & 3   &  11 & (0.2224, 0.0159) & (1.7110,  0.3100) &  0.1117 &  0.2067 \\
Holway  & 3   &   9 & (0.2169, 0.0215) & (1.6974,  0.2944) &  0.0189 &  0.1126 \\ 
Holway*  & 3   &   9 & (0.2149, 0.0225) & (1.6927,  0.2901) &  0.0168 &  0.0924 \\ 
 MM & 3   & 201 & (0.2067, -0.0004) & (1.6844,  0.2757) &  0.0000 &  0.0939 \\ 
\hline
\end{tabular}
\caption{Summary of model quality for the improved exponential model when applied to data from Morozov et al.~\cite{MGG-24} (tests 2--5) and Frezzotti~\cite{Fr-07} (tests 7--9, 11--14).}
\label{tab:fitGen}
\end{table}
\begin{figure}
    \centering
    \includegraphics[scale=1]{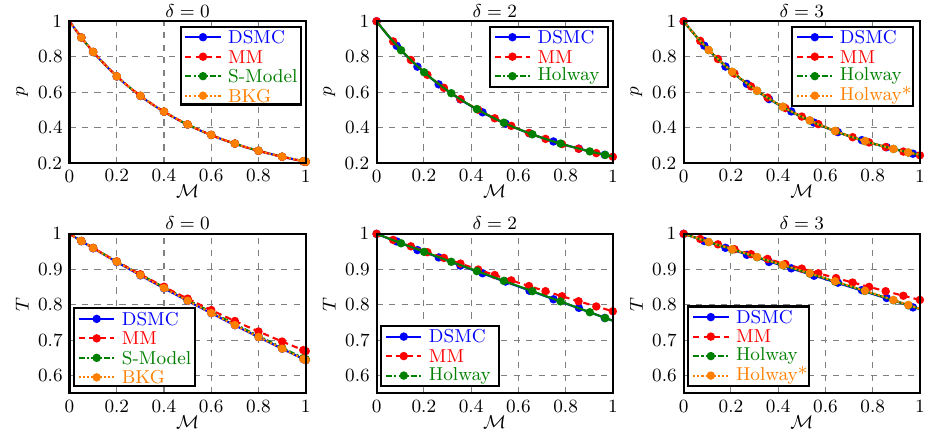}
    \caption{Fitted models to data from~\cite{MGG-24} ($\delta=0$) and~\cite{Fr-07} ($\delta\in\{2,3\}$).}
    \label{fig:fitGen}
\end{figure}

The data in Table~\ref{tab:fitGen} show the MAPE values obtained when fitting the data of Morozov et al.~\cite{MGG-24} and Frezzotti~\cite{Fr-07} using the improved exponential model~\eqref{eq:fit2}, together with the corresponding values obtained in this work, considering the maximal total entropy production (MTEP).
The data sets considered by Frezzotti~\cite{Fr-07} and Morozov et al.~\cite{MGG-24} were obtained by applying the Direct Simulation Monte Carlo (DSMC) method, considering the Bhatnagar–Gross–Krook (BGK) operator; here the data sets are originally from Sone et al.~\cite{SS-90,So-00}, the S-model, moment methods (MM), or Holway's kinetic model. 
For polyatomic gases, $\delta=2$ or $\delta=3$, Frezzotti~\cite{Fr-07} considered collision operators where a fraction of $0.3$ of the collisions was inelastic and, correspondingly, a fraction of $0.7$ of the collisions was elastic; the only exception is the data set based on the Holway model (*) for $\delta=3$ where all collisions are assumed to be inelastic~\cite{Fr-07}. 
However, note that the results obtained in the present work (MTEP) do not depend on the specific collision operator---a fraction of the collisions, but not all of them, can be resonant.
Figure~\ref{fig:fitGen} shows fitted curves together with the data, represented by circles at all points in the small data sets $(n \in \left\{9,11,12\right\})$ and a subset of points in the large data set ($n=201$).

\begin{figure}
    \centering
    \includegraphics[scale=1]{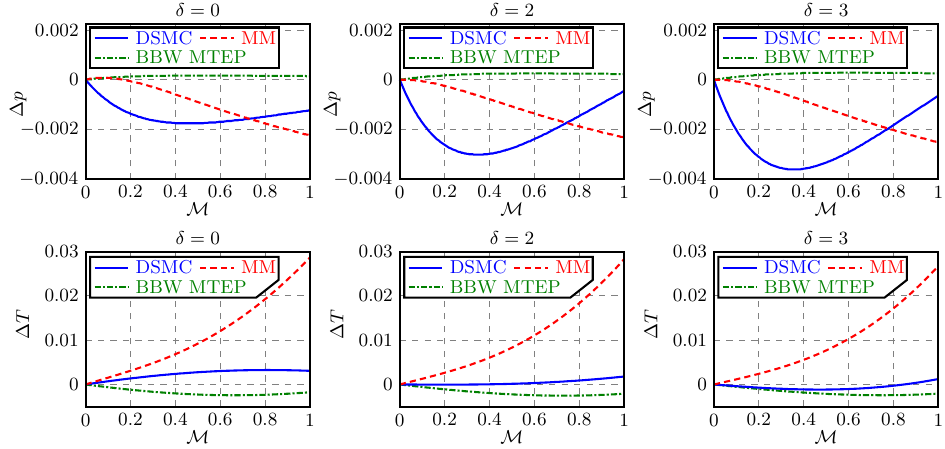}
    \caption{Differences of some fitted models---for data obtained by the DSMC method, moment methods (MM), and the maximal total entropy production (MTEP) considered in~\cite{BBW-26}---and the fitted model for the modified maximal total entropy production considered in this paper.}
    \label{fig:fitComp}
\end{figure}
The differences of the fitted models for the data obtained based on the DSMC and moment methods from Frezzotti~\cite{Fr-07} and Morozov et al.~\cite{MGG-24}, as well as for the maximal total entropy production (MTEP) considered in our article~\cite{BBW-26}, and the fitted model for the modified maximal total entropy production considered in this article are shown in Figure~\ref{fig:fitComp}. The zero line corresponds to perfect correspondence with the fitted model for the modified MTEP. Note that the absolute difference never exceeds 0.004 for the models based on the DSMC-method, neither for $p$ nor for $T$.

\section{Generalized exponential model fitted to the ratio of specific heats \label{GEM}}
\begin{figure}
\centering
\includegraphics[scale=1]{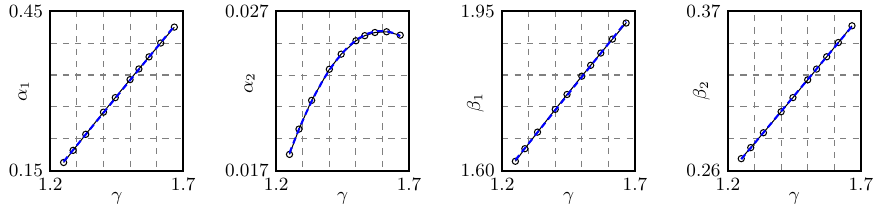}
\caption{Optimized parameters for the improved exponential model~\eqref{eq:fit2} as functions of the ratio of specific heats $\gamma$ for $\delta\in\{0, 0.25, 0.5, 0.75, 1, 1.5, 2, 3, 4, 5\}$ marked by circles connected by thin solid lines together with least-squares fitted trend functions; linear functions for $\alpha^\ddag_1$, $\beta^\ddag_1$, and $\beta^\ddag_2$, and a quadratic function for $\alpha^\ddag_2$ drawn by thick dashed lines.}
\label{fig:fit2}
\end{figure}
The optimized parameters in Table~\ref{tab:fit2} suggest a systematic dependence on $\delta$ for $\alpha^\delta_{\ddag}$ and $\beta^\delta_\ddag$.
Although $\delta$ is a natural parameter from some perspectives, several thermodynamic properties are expressed more conveniently in terms of the ratio of specific heats $\gamma$~\eqref{rsh}.
Expressing the parameter dependence in terms of $\gamma$ links the model to a physically significant quantity and reveals remarkably simple trends in the fitted coefficients.

Figure~\ref{fig:fit2} illustrates this dependence by showing the optimized parameters of models~\eqref{eq:fit2} as functions of $\gamma$ when using data for $\delta\in\{0, 0.25, 0.5, 0.75, 1, 1.5, 2, 3, 4, 5\}$.
In particular, the observed trends indicate that $\alpha_1$, $\beta_1$, and $\beta_2$ depend approximately linearly on $\gamma$, while $\alpha_2$ exhibits a noticeable nonlinear  
(apparently quadratic) dependence. Moreover, the magnitude and range of $\alpha_2$ are significantly smaller than those of $\alpha_1$, $\beta_1$, and $\beta_2$.

\begin{figure}
\centering
\includegraphics[scale=1]{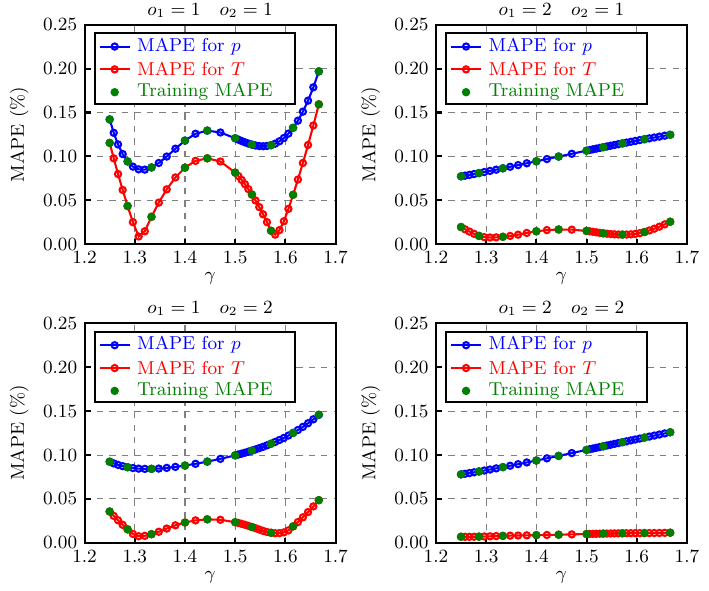}
\caption{
Mean absolute percentage error (MAPE) as functions of $\gamma$ for the $\gamma$-dependent model~\eqref{eq:fit3} with polynomial coefficients~\eqref{eq:polynomialCoeffs}, fitted to maximal total entropy production data for $\delta\in\{0, 0.25, 0.5, 0.75, 1, 1.5, 2, 3, 4, 5\}$ (training data). 
Each panel corresponds to a choice of polynomial degrees $(o_1,o_2)\in\{1,2\}\times\{1,2\}$, where $o_1:=o_1^\alpha=o_1^\beta$ and $o_2:=o_2^\alpha=o_2^\beta$.
}
\label{fig:fit3MAPEOrders}
\end{figure}

To enable prediction for arbitrary degrees of freedom, we construct a generalized model in which the parameters of the exponential ansatz are represented as smooth functions of the ratio of specific heats $\gamma$.
Motivated by the trends discussed above, we model the coefficients in model~\eqref{eq:fit2} as low-order polynomials in $\gamma$.
We thus consider the ansatz
\begin{equation}
  p_\ast(\mathcal{M},\gamma) = \exp\!\Big( -\mathcal{B}_1(\gamma)\,\mathcal{M} + \mathcal{B}_2(\gamma)\,\mathcal{M}^2 \Big)
  \quad\text{and}\quad
  T_\ast(\mathcal{M},\gamma) = \exp\!\Big( -\mathcal{A}_1(\gamma)\,\mathcal{M} - \mathcal{A}_2(\gamma)\,\mathcal{M}^3 \Big),
  \label{eq:fit3}
\end{equation}
where, for $i\in \{1,2\}$,
\begin{equation}
    \mathcal{A}_i(\gamma) = \sum_{j=0}^{o^\alpha_i} \alpha_{ij}\gamma^j 
    \quad\text{and}\quad
    \mathcal{B}_i(\gamma) = \sum_{j=0}^{o^\beta_i} \beta_{ij}\gamma^j,
    \label{eq:polynomialCoeffs}
\end{equation}
with polynomial degrees $o^\alpha_i$ and $o^\beta_i$ to be determined.
As a baseline model, we use $o^\alpha_i = o^\beta_i =1$, for $i\in\{1,2\}$.
To determine suitable polynomial degrees in~\eqref{eq:polynomialCoeffs}, we evaluate the performance of the generalized model~\eqref{eq:fit3} for different combinations of polynomial degrees. 
For each order combination, the coefficient functions are determined by fitting the generalized model to the maximal total entropy production curves for $\delta\in\{0, 0.25, 0.5, 0.75, 1, 1.5, 2, 3, 4, 5\}$ sampled at $\mathcal{M}\in\{0.01,0.02,\ldots,1\}$ in the least-squares sense.
Figure~\ref{fig:fit3MAPEOrders} shows the resulting mean absolute percentage errors (MAPE) as functions of $\gamma$ for both the relative pressure $p$ and the relative temperature $T$.

These results show that increasing the polynomial degrees associated with  $\mathcal{A}_1$ and $\mathcal{B}_1$ provides a significant improvement in accuracy compared to the baseline model, while increasing the degree associated with $\mathcal{A}_2$ and $\mathcal{B}_2$ provides a more moderate reduction in error. 
In particular, for the relative pressure $p$, the accuracy has essentially saturated once $\mathcal{B}_1$ is quadratic, and allowing $\mathcal{B}_2$ to be quadratic, or even of higher order, in $\gamma$ does not lead to a significant reduction in the MAPE for $p$.
However, for the relative temperature $T$, a small but consistent improvement is observed when increasing the degree associated with  $\mathcal{A}_2$, even after increasing the degree associated with  $\mathcal{A}_1$.
Thus, unlike $\mathcal{B}_2$, using a quadratic dependence on $\gamma$ for $\mathcal{A}_2$ provides a measurable improvement in the approximation of $T$.
Tests with higher polynomial orders confirm that the MAPE has saturated and that no further improvement is obtained by increasing the degrees beyond those considered here.

\begin{figure}
\centering
\includegraphics[scale=1]{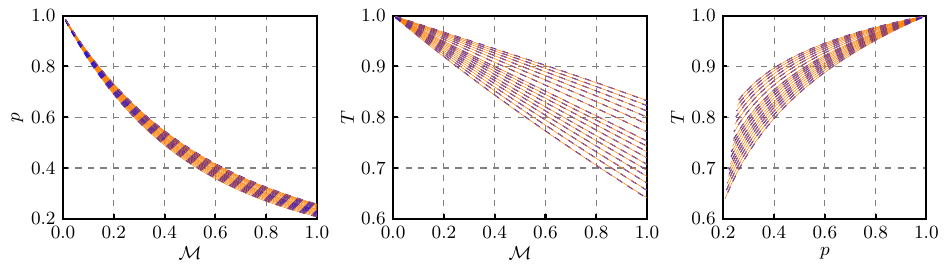}
\caption{The maximal total entropy production curves (dashed lines) and the corresponding predictions of the generalized\\ 
$\gamma$-dependent model~\eqref{eq:fit3} (solid lines) obtained using the optimized coefficient functions in expression~\eqref{eq:fit3OptCoeffs} for the values\\ $\delta\in\{0, 0.2, 0.4, 0.6, 0.8, 1, 1.25, 1.5, 1.75, 2, 2.5, 3, 3.5, 4, 4.5, 5\}$.
The larger the value of $\delta$, the higher the corresponding curve.
}
\label{fig:fit3MpT}
\end{figure}

Table~\ref{tab:fit3MAPE} complements Figure~\ref{fig:fit3MAPEOrders} by reporting the average MAPE values for selected combinations of polynomial degrees, including the case $o_1^\alpha=o_2^\alpha=o_1^\beta=o_2^\beta=3$, which illustrates the saturation of the model accuracy at low polynomial orders.
Based on these observations, we select polynomial degrees corresponding to the onset of saturation, taking $o^\alpha_1 = o^\alpha_2 = o^\beta_1 = 2$ and $o^\beta_2 = 1$, corresponding to quadratic polynomials for $\mathcal{A}_1$, $\mathcal{A}_2$, and $\mathcal{B}_1$, and a linear polynomial for $\mathcal{B}_2$.
With these selected polynomial degrees, the optimized coefficient functions in the least-squares sense are as follows
\begin{equation}
\begin{aligned}
\mathcal{A}_1^\ast(\gamma) &= -0.7529 + 0.8272 \gamma -0.0740 \gamma^2, &
\mathcal{A}_2^\ast(\gamma) &= -0.1386 + 0.2062 \gamma -0.0647 \gamma^2, \\
\mathcal{B}_1^\ast(\gamma) &= 0.4516 +  1.0905 \gamma -0.1243 \gamma^2, &
\mathcal{B}_2^\ast(\gamma) &= -0.0068 + 0.2198 \gamma. \\
\end{aligned}
\label{eq:fit3OptCoeffs}
\end{equation}
These optimized coefficient functions define a single $\gamma$-dependent model valid for all considered ratios of specific heats.
For the selected model, the average MAPE value is below $0.01\,\%$ for $T$ and approximately $0.10\,\%$ for $p$, demonstrating that the generalized model retains excellent accuracy while substantially increasing its applicability range.
Figure~\ref{fig:fit3MpT} compares the maximal total entropy production curves (dashed lines) for $\delta\in\{0, 0.2, 0.4, 0.6, 0.8, 1, 1.25, 1.5, 1.75, 2, 2.5, 3, 3.5, 4, 4.5, 5\}$
with the corresponding predictions of the generalized $\gamma$-dependent model~\eqref{eq:fit3} (solid lines) obtained using the optimized coefficient functions in expression~\eqref{eq:fit3OptCoeffs}.
The agreement also remains excellent for intermediate values of $\delta$ not included in the training data, indicating that the generalized model provides accurate interpolation throughout the range of specific heat ratios considered.
\begin{table}
\centering
\begin{tabular}{c|ccccc}
 $(o_1^\alpha,o_2^\alpha)=(o_1^\beta,o_2^\beta)$ & (1,1) & (1,2) & (2,1) & (2,2) & (3,3) \\
\hline
MAPE for $T$ & 0.0634 & 0.0204 & 0.0136 & 0.0096 & 0.0095 \\
MAPE for $p$ & 0.1197 & 0.1040 & 0.1030 & 0.1030 & 0.1029 \\
\end{tabular}
\caption{Average MAPE values for the $\gamma$-dependent model~\eqref{eq:fit3} with polynomial coefficients~\eqref{eq:polynomialCoeffs}, fitted to maximal total entropy production data for $\delta\in\{0, 0.25, 0.5, 0.75, 1, 1.5, 2, 3, 4, 5\}$ (training data) for selected combinations of polynomial degrees.}
\label{tab:fit3MAPE}
\end{table}

\section{Mixed diffuse reflection and complete absorption boundary conditions \label%
{SGBC}}
In this section, we consider mixed diffuse reflection and complete absorption boundary conditions~\cite{STG-01,Fr-07}.
Due to mass conservation, we obtain by applying the boundary conditions~\eqref{BC1} that
\begin{equation*}
\int_{\mathbb{R}_{-}^{3}}\xi _{1}f\,d\boldsymbol{\xi }=n_{\infty }u-\sqrt{%
\frac{k_{\mathrm{B}}T_{0}}{2\pi m}}\left( \sigma _{e}n_{0}+(1-\sigma
_{e})N_{0}\right) \text{.}
\end{equation*}%
Then%
\begin{equation*}
N_{0}=-\sqrt{\frac{2\pi m}{k_{\mathrm{B}}T_{0}}}\int_{\mathbb{R}_{-}^{3}}\xi
_{1}f\,\mathrm{d}\boldsymbol{\xi }=\left( \sigma _{e}n_{0}+(1-\sigma
_{e})N_{0}\right) -\sqrt{\frac{2\pi m}{k_{\mathrm{B}}T_{0}}}n_{\infty }u%
\text{,}
\end{equation*}%
and solving for $N_{0}$, we obtain the expression 
\begin{equation*}
N_{0}=n_{0}-\sqrt{\frac{2\pi m}{k_{\mathrm{B}}T_{0}}}\frac{n_{\infty }}{%
\sigma _{e}}u\text{ for }0<\sigma _{e}\leq 1\text{.}
\end{equation*}%
Hence,%
\begin{equation*}
\frac{\sigma _{e}n_{0}+(1-\sigma _{e})N_{0}}{n_{0}}=1-\dfrac{1-\sigma _{e}}{%
\sigma _{e}}\sqrt{\dfrac{2\pi m}{k_{\mathrm{B}}T_{0}}}\frac{n_{\infty }}{%
n_{0}}u=1-\dfrac{1-\sigma _{e}}{\sigma _{e}}\sqrt{\frac{2\pi \gamma }{T}}p%
\mathcal{M}\text{.}
\end{equation*}%
Denote by%
\begin{equation*}
\widetilde{p}_{0}=\left( \sigma _{e}n_{0}+(1-\sigma _{e})N_{0}\right) k_{%
\mathrm{B}}T_{0}\text{ and }\widetilde{p}=\frac{p_{\infty }}{\widetilde{p}%
_{0}}\text{.}
\end{equation*}%
Then%
\begin{equation*}
\frac{\widetilde{p}_{0}}{p_{0}}=  \frac{\sigma _{e}n_{0}+(1-\sigma
_{e})N_{0}}{n_{0}} =1-\dfrac{1-\sigma _{e}}{\sigma _{e}}\sqrt{\frac{2\pi
\gamma }{T}}\widetilde{p}\mathcal{M}\frac{\widetilde{p}_{0}}{p_{0}}
\end{equation*}%
and solving for $p_{0}/\widetilde{p}_{0}$, we obtain the expression 
\begin{equation*}
\frac{p_{0}}{\widetilde{p}_{0}}=1+\dfrac{1-\sigma _{e}}{\sigma _{e}}\sqrt{%
\frac{2\pi \gamma }{T}}\widetilde{p}\mathcal{M}>0\text{,}
\end{equation*}%
and we can conclude that%
\begin{equation*}
p=\widetilde{p}\frac{\widetilde{p}_{0}}{p_{0}}=\widetilde{p}\left( 1+\dfrac{%
1-\sigma _{e}}{\sigma _{e}}\sqrt{\frac{2\pi \gamma }{T}}\widetilde{p}%
\mathcal{M}\right) ^{-1}\text{.}
\end{equation*}

By Proposition~\ref{P1}, we obtain the following proposition.

\begin{figure}
\centering
\includegraphics[scale=1]{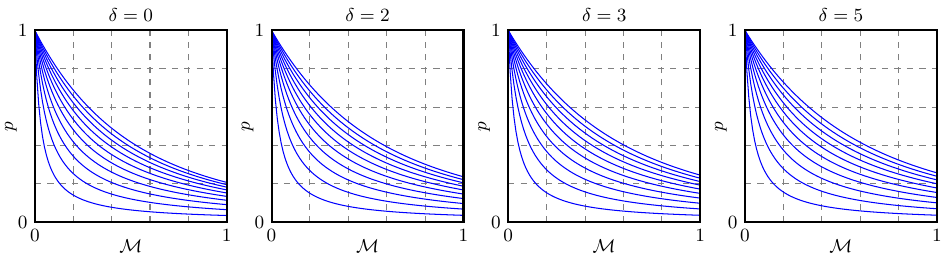}
\caption{Value of $p$ that maximizes the total entropy production as a function of the Mach number $\mathcal{M}$ for different values of $\sigma _{e}\in\left\{0.1, 0.2, \ldots, 1.0\right\}$---the larger the value of $\sigma _{e}$, the higher the corresponding curve---and $\delta\in\left\{0,2,3,5\right\}$.}
\label{fig:alphaScaling}
\end{figure}

\begin{proposition}
For the half-space problem~\eqref{HSP1}, with boundary conditions~\eqref{BC1}, to admit a solution, the following relations for the parameters of the Maxwellians at the condensed interface and the uniform phase at infinity have to be fulfilled:

i) for all $u>0$ or equivalently for all $\mathcal{M}_{\infty }>0$
\begin{align*}
p\mathcal{M} &\leq \sigma _{e}\sqrt{\frac{T}{2\pi \gamma }}\text{,}  \notag \\
p\mathcal{M}\left( \left( 3+\delta +\mathcal{M}^{2}\right) T+%
\dfrac{1-\sigma _{e}}{\sigma _{e}}\frac{4+\delta }{\gamma }\right)  &\leq %
\frac{4+\delta }{\gamma }\sqrt{\frac{T}{2\pi \gamma }}\text{,}  \notag \\
\frac{1}{2\left( 1+\gamma \mathcal{M}^{2}\right) +\dfrac{1-\sigma _{e}}{%
\sigma _{e}}\sqrt{\dfrac{2\pi \gamma }{T}}\mathcal{M}} \leq
p&\leq \frac{1}{\left( 1+\dfrac{\mathcal{M}^{2}}{3+\delta }\right) ^{\left(
5+\delta \right) /2}+\dfrac{1-\sigma _{e}}{\sigma _{e}}\sqrt{\dfrac{2\pi
\gamma }{T}}\mathcal{M}}\text{;}  \label{NCP1}
\end{align*}

ii) for $u=0$ or equivalently $\mathcal{M}_{\infty }=0$
\begin{equation}
p=T=1\text{.}  \label{TC1}
\end{equation}%
Moreover, under condition~\eqref{TC1} there exists a unique
solution
\begin{equation*}
f=f(x,\boldsymbol{\xi },I)=M_{0}(\boldsymbol{\xi },I)\text{.} 
\end{equation*}%
\end{proposition} 
The value of $p$ that maximizes the total entropy production for different values of $\sigma _{e}$ is shown in Figure~\ref{fig:alphaScaling} for $\delta=0$, $2$, $3$, and $5$, while the corresponding value of $T$ is independent of $\sigma_e$ and coincides with that of the complete absorption case. 
We observe that the difference between the curves corresponding to different internal degrees of freedom decreases as $\sigma _{e}\rightarrow 0$.

\section{Concluding Discussion\label{CD}}
In this work, we revisited the half-space problem of steady evaporation for polyatomic gases using a modified total entropy production functional.
By analyzing how the evaporation curves depend on the Mach number and the ratio of specific heats, we identified simple structures in the pressure and temperature ratios and used these observations to construct compact analytical approximations.

The principal finding is the remarkably simple structure underlying these maximal total entropy production curves.
Starting from exponential representations in the Mach number of the pressure and temperature curves corresponding to maximal total entropy production for fixed numbers of internal degrees of freedom, we found that the remaining discrepancies are highly structured and can be accounted for using only a few additional parameters.
Moreover, the resulting coefficients vary smoothly with the ratio of specific heats, making it possible to construct a generalized $\gamma$-dependent model that accurately represents the corresponding maximal total entropy production curves.
The resulting model achieves this accuracy while retaining a compact form, suggesting that the pressure and temperature curves corresponding to maximal total entropy production exhibit a surprisingly low-dimensional dependence on the Mach number and the ratio of specific heats.

The proposed functional forms were also found to accurately represent the evaporation data obtained by several alternative numerical and kinetic approaches. 
Furthermore, the extension to mixed diffuse reflection and complete absorption boundary conditions shows that the proposed modeling approach remains applicable beyond the setting primarily considered in this work. 
Together, these results indicate that the observed structure may capture broader features of steady evaporation phenomena.

\section*{Acknowledgement}
The authors thank Professor A.~Frezzotti for sharing the data sets from his study~\cite{Fr-07}.

\end{document}